\documentclass[a4paper,12pt]{article}

\usepackage[english]{babel} 
\usepackage[ansinew]{inputenc}
\usepackage[LGR,T1]{fontenc}

\usepackage{verbatim}
\usepackage{color}
\usepackage{calc}
\usepackage{hyperref}
\usepackage[dvips]{graphicx}
\hypersetup{colorlinks=true, linkcolor=blue, filecolor=blue, urlcolor=blue}

\usepackage{amssymb}
\usepackage[mathscr]{eucal}
\usepackage{amscd,amsmath,amsfonts,amssymb,textcomp}
\everymath{\displaystyle}

\usepackage[all]{xy}

\newcommand{\MBC}{\mathbb{C}}

\newcommand{\MBN}{\mathbb{N}}

\newcommand{\MBQ}{\mathbb{Q}}

\newcommand{\MBZ}{\mathbb{Z}}

\newcommand{\MCC}{\mathcal{C}}
\newcommand{\MCD}{\mathcal{D}}
\newcommand{\MCE}{\mathcal{E}}

\newcommand{\MCG}{\mathcal{G}}
\newcommand{\MCH}{\mathcal{H}}

\newcommand{\MCK}{\mathcal{K}}

\newcommand{\MCO}{\mathcal{O}}

\newcommand{\MCS}{\mathcal{S}}

\newcommand{\MCU}{\mathcal{U}}

\newcommand{\MCZ}{\mathcal{Z}}

\newcommand{\MFa}{\mathfrak{a}}

\newcommand{\MFf}{\mathfrak{f}}

\newcommand{\MFg}{\mathfrak{g}}

\newcommand{\MFm}{\mathfrak{m}}

\newcommand{\MFp}{\mathfrak{p}}

\newcommand{\MFq}{\mathfrak{q}}

\newcommand{\MFr}{\mathfrak{r}}

\newcommand{\MFu}{\mathfrak{u}}

\newcommand{\MSN}{\mathscr{N}}

\newcommand{\SCm}{\textsc{m}}

\newcommand{\SCr}{\textsc{r}}

\newcommand{\GGa}{\alpha}

\newcommand{\GGd}{\delta}

\newcommand{\GVf}{\varphi}
\newcommand{\GGg}{\gamma}
\newcommand{\GGG}{\Gamma}

\newcommand{\GGl}{\lambda}
\newcommand{\GGL}{\Lambda}

\newcommand{\GGs}{\sigma}

\newcommand{\GGz}{\zeta}

\newcommand{\ba}[2]{{\mathop{#1}\limits_{#2}}}

\newcommand{\lA}{\left\{}

\newcommand{\alg}{\mathrm{alg}}
\newcommand{\Char}{\mathrm{char}}

\newcommand{\Cl}{\mathrm{Cl}}

\newcommand{\Cok}{\mathrm{Cok}}

\newcommand{\Gal}{\mathrm{Gal}}

\newcommand{\IM}{\mathrm{Im}}
\newcommand{\Ker}{\mathrm{Ker}}

\newcommand{\tor}{\mathrm{tor}}

\newcommand{\II}{\mathrm{I}}
\newcommand{\MM}{\mathrm{M}}

\usepackage{makeidx}
\title{GLOBAL UNITS MODULO ELLIPTIC UNITS AND $2$-IDEAL CLASS GROUPS}
\author{\textsc{St\'ephane VIGUI\'E}
\footnote{St\'ephane Viguié, Laboratoire de mathématiques de Besançon, UMR CNRS 6623, Université de Franche-Comté, 16 route de Gray, 25030 Besançon cedex, France.
e-mail: \texttt{stephane.viguie@univ-fcomte.fr}}
}
\makeindex

\newtheorem{cro}{cro}[section]
\newtheorem{dfe}{dfe}[section]
\newtheorem{lem}{lem}[section]
\newtheorem{pro}{pro}[section]
\newtheorem{rem}{rem}[section]
\newtheorem{teh}{teh}[section]

\newtheorem{cor}[cro]{Corollary}

\newtheorem{df}[dfe]{Definition}
\newtheorem{lm}[lem]{Lemma}
\newtheorem{pr}[pro]{Proposition}

\newtheorem{rmq}[rem]{Remark}
\newtheorem{theo}[teh]{Theorem}

\numberwithin{equation}{section}

\begin{document}
\maketitle


\begin{abstract}
Let $p\in\{2,3\}$, and let $k$ be an imaginary quadratic number field in which $p$ decomposes into two distinct primes $\MFp$ and $\bar{\MFp}$.
Let $k_\infty$ be the unique $\MBZ_p$-extension of $k$ which is unramified outside of $\MFp$, and let $K_\infty$ be a finite extension of $k_\infty$, abelian over $k$.
We prove that in $K_\infty$, the projective limit of the $p$-class group and the projective limit of units modulo elliptic units share the same $\mu$-invariant and the same $\GGl$-invariant.
Then we prove that up to a constant, the characteristic ideal of the projective limit of the $p$-class group coincides with the characteristic ideal of the projective limit of units modulo elliptic units. 
\end{abstract}
\bigskip

\\


\section{Introduction.}

Let $p$ be a prime number, and let $k$ be an imaginary quadratic number field in which $p$ decomposes into two distinct primes $\MFp$ and $\bar{\MFp}$.
Let $k_\infty$ be the unique $\MBZ_p$-extension of $k$ which is unramified outside of $\MFp$, and let $K_\infty$ be a finite extension of $k_\infty$, abelian over $k$.
Let $G_\infty$ be the Galois group of $K_\infty/k$.
We choose a decomposition of $G_\infty$ as a direct sum of a finite group $G$ (the torsion subgroup of $G_\infty$) and a topological group $\GGG$ isomorphic to $\MBZ_p$, $G_\infty = G\times \GGG$.
For any $n\in\MBN$, let $K_n$ be the field fixed by $\GGG_n := \GGG^{p^n}$, and let $G_n:=\Gal\left(K_n/k\right)$.
Remark that there may be different choices for $\GGG$, but when $p^n$ is larger than the order of the $p$-part of $G$, the group $\GGG_n$ does not depend on the choice of $\GGG$.

Let $F/k$ be an abelian extension of $k$.
If $[F:k]<\infty$, we denote by $\MCO_F$ the ring of integers of $F$.
Then we write $\MCO_F^\times$ for the group of global units of $F$, and $C_F$ for the group of elliptic units of $F$ (see section \ref{sectionellipticunits}).
Also we let $A_F$ be the $p$-part of the class group $\Cl\left(\MCO_F\right)$ of $\MCO_F$.
We set $\MCE_F:=\MBZ_p\otimes_\MBZ\MCO_F^\times$ and $\MCC_F:=\MBZ_p\otimes_\MBZ C_F$.
When $F/k$ is infinite, we define $\MCE_F$, $\MCC_F$ and $A_F$, by taking projective limits over finite sub-extensions, under the norm maps.
For any $n\in\MBN\cup\{\infty\}$, we set $\MCE_n:=\MCE_{K_n}$, $\MCC_n:=\MCC_{K_n}$, and $A_n:=A_{K_n}$.

For any profinite group $\MCG$, and any commutative ring $R$, we define the Iwasawa algebra
\[R\left[\left[\MCG\right]\right] := \varprojlim R\left[\MCH\right],\]
where the projective limit is over all finite quotient $\MCH$ of $\MCG$.
In case $\MCG=G_\infty$ (respectively $\MCG=\GGG$), we shall write
\[\GGL := \MBZ_p\left[\left[G_\infty\right]\right] \quad \left(\text{respectively} \quad \GGL' := \MBZ_p\left[\left[\GGG\right]\right] \right).\]
Then $A_\infty$ and $\MCE_\infty/\MCC_\infty$ are naturally $\GGL$-modules.
As we shall see below, they are finitely generated and torsion over $\GGL'$.
Let us fix a topological generator $\GGg$ of $\GGG$, and set $T:=\GGg-1$.
Then for any finite extension $L/\MBQ_p$, $\MCO_L[[\GGG]]$ is isomorphic to $\MCO_L[[T]]$, where $\MCO_L$ is the ring of integers of $L$.
It is well known that $\MCO_L[[T]]$ is a noetherian, regular, local domain.
We also recall that $\MCO_L[[T]]$ is a unique factorization domain.
If $u$ is a uniformizer of $\MCO_L$, then the maximal ideal $\MFm$ of $\MCO_L$ is generated by $u$ and $T$, and $\MCO_L[[T]]$ is a compact topological ring with respect to its $\MFm$-adic topology.
A morphism $f:M\rightarrow N$ between two finitely generated $\MCO_L[[T]]$-module is called a pseudo-isomorphism if its kernel and its cokernel are finite.
If a finitely generated $\MCO_L[[T]]$-module $M$ is given, then one may find elements $P_1$, ..., $P_r$ in $\MCO_L[T]$, irreducible in $\MCO_L[[T]]$, and nonnegative integers $n_0$, ..., $n_r$, such that there is a pseudo-isomorphism
\[M\longrightarrow \MCO_L[[T]]^{n_0} \oplus \bigoplus_{i=1}^r \MCO_L[[T]]/\left(P_i^{n_i}\right).\]
Moreover, the integer $n_0$ and the ideals $\left(P_1^{n_1}\right)$, ..., $\left(P_r^{n_r}\right)$, are uniquely determined by $M$.
If $n_0=0$, then the ideal generated by $P_1^{n_1}\cdots P_r^{n_r}$ is called the characteristic ideal of $M$, and is denoted by $\Char_{\MCO_L[[T]]}(M)$.

We denote by $\MBC_p$ a completion of an algebraic closure of $\MBQ_p$.
Let $\chi:G\rightarrow\MBC_p$ be an irreducible character of $G$.
Let $\MBQ_p(\chi)\subset\MBC_p$ be the abelian extension of $\MBQ_p$ generated by the values of $\chi$.
We denote by $\MBZ_p(\chi)$ the ring of integers of $\MBQ_p(\chi)$.
The group $G$ acts naturally on $\MBQ_p(\chi)$ if we set, for all $g\in G$ and all $x\in\MBQ_p(\chi)$, $g.x:=\chi(g)x$.
For any $\MBZ_p[G]$-module $Y$, we define the $\chi$-quotient $Y_\chi$ by $Y_\chi:=\MBZ_p(\chi)\otimes_{\MBZ_p[G]}Y$.
If $Y$ is a $\GGL$-module, then $Y_\chi$ is a $\MBZ_p(\chi)[[T]]$-module in a natural way.
As a particular case, $\GGL_\chi \simeq \MBZ_p(\chi)[[T]]$.
For any finitely generated $\GGL_\chi$-module $Z$, we shall denote $\Char_{\GGL_\chi}Z$ simply by $\Char_\chi Z$.

The goal of this article is to prove Theorem \ref{conjmain} below, which is a raw formulation of the (one-variable) main conjecture at $p=2$ or $p=3$.
In \cite[Theorem 4.1]{rubin91} and \cite[Theorem 2]{rubin94}, Rubin used Euler systems to prove the main conjectures for $\MBZ_p$ or $\MBZ_p^2$ extensions of a finite abelian extension $F$ of $k$, where $p\nmid w_k[F:k]$, $w_k$ being the number of roots of unity in $k$.
Inspired by the ideas of Rubin, Greither used Euler systems to prove the main conjecture for cyclotomic units and for the cyclotomic $\MBZ_p$-extension $F_\infty/F$, with $F_\infty$ abelian over $\MBQ$ (see \cite[Theorem 3.2]{greither92}).
Bley proved the main conjecture when $p\nmid2\#\left(\Cl\left(\MCO_k\right)\right)$, and when there is a nonzero ideal $\MFf$ of $\MCO_k$, prime to $\MFp$, 
such that for all $n\in\MBN$, $K_n = k\left(\MFf\MFp^n\right)$ is the ray class field of $k$ modulo $\MFf\MFp^n$ (see \cite[Theorem 3.1]{bley06}).
More recently, Hassan Oukhaba adapted Rubin's method and obtained the divisibility relation 
\begin{equation}
\label{divrel}
\Char_\chi \left( A_\infty \right) \vert p^{m_\chi}\Char_\chi \left( \MCE_\infty / \MCC_\infty \right) \quad\text{for some $m_\chi\in\MBN$,}
\end{equation}
for $p=2$, still under the condition $2\nmid\left[K_0:k\right]$ (see \cite{oukhaba10}).
In \cite{viguie11b}, we proved the main conjecture for the extension $K_\infty / K$ when $p\notin\{2,3\}$.
When $p\in\{2,3\}$, we obtained the divisibility relation (\ref{divrel}).
Here we prove that an equality holds also for $p\in\{2,3\}$, up to a constant in $\MBZ_p(\chi)$.

\begin{theo}\label{conjmain}
Let $p\in\{2,3\}$, and let $\chi$ be an irreducible $\MBC_p$ character on $G$.
There are $a_\chi\in\MBN$ and $b_\chi\in\MBN$ such that 
\begin{equation}
\label{divisibility}
\MFu_\chi^{a_\chi} \Char_\chi\left(A_{\infty,\chi}\right) = \MFu_\chi^{b_\chi} \Char_\chi\left(\MCE_\infty/\MCC_\infty\right)_\chi,
\end{equation}
where $\MFu_\chi$ is a uniformizer of $\MBZ_p(\chi)$.
\end{theo}

Let $L/\MBQ_p$ be a finite algebraic extension.
For any finitely generated torsion $\MCO_L[[T]]$-module $M$, we recall that the $\GGl$-invariant $\GGl(M)$ of $M$ is the $\MCO_L$-rank of $M$, which is also the Weierstrass degree of any generator of $\Char_{\MCO_L[[T]]}M$.
Let $\MFu_L$ be a uniformizer of $\MCO_L$.
Then the $\mu$-invariant $\mu(M)$ of $M$ is the maximal power of $\MFu_L$ which divides $\Char_{\MCO_L[[T]]}M$.

For any finitely generated $\GGL$-module $M$, we recall that
\begin{equation}
\label{somlambdachi}
\sum_{\chi}\,\GGl\left(M_\chi\right) \; = \; \GGl\left(M\right),
\end{equation}
where the sum is over all irreducible $\MBC_p$ characters of $G$.
In particular, (\ref{somlambdachi}) and the divisibility relation (\ref{divrel}) gives us
\begin{equation}
\label{OuGr}
\GGl\left(A_\infty\right) \; = \; \sum_{\chi}\,\GGl\left(A_{\infty,\chi}\right) \; \leq \; \sum_{\chi}\,\GGl\left(\MCE_{\infty,\chi} / \MCC_{\infty,\chi}\right) \; = \; \GGl\left(\MCE_{\infty} / \MCC_{\infty}\right),
\end{equation}
where the sums are over all irreducible $\MBC_p$ characters of $G$.
Hence in order to prove Theorem \ref{conjmain}, it suffices to show that 
\begin{equation}
\GGl \left( \MCE_\infty / \MCC_\infty \right) \, = \, \GGl \left( A_\infty \right).
\label{MKjbvVT}
\end{equation}
Indeed from (\ref{MKjbvVT}) and (\ref{OuGr}) we deduce $\GGl\left(A_{\infty,\chi}\right) = \GGl\left(\MCE_{\infty,\chi} / \MCC_{\infty,\chi}\right)$ for any $\MBC_p$ character $\chi$ of $G$, and this last equality together with (\ref{divrel}) induce Theorem \ref{conjmain}.
For $p\notin\{2,3\}$, a proof of (\ref{MKjbvVT}) is sketched in \cite[III.2.1]{deshalit87}.
A complete proof is given in the cyclotomic case by J.\,R.\,Belliard in \cite{belliard09}.
In this article, we prove the elliptic version of \cite[Theorem 6.2]{belliard09}.
More precisely we show the two following facts (see Theorem \ref{theoegalGGlmu}).

$\mathrm{(i)}$ The $\GGL'$-modules $\MCE_\infty / \MCC_\infty$ and $A_\infty$ share the same $\GGl$-invariant and the same $\mu$-invariant.

$\mathrm{(ii)}$ The $\GGL'$-modules $\MCU_\infty / \MCC_\infty$ and $B_\infty$ share the same $\GGl$-invariant and the same $\mu$-invariant, 
where $\MCU_\infty$ is a module of semi-local units and $B_\infty$ is the Galois group of the maximal abelian pro-$p$-extension of $K_\infty$ which is unramified outside the primes above $\MFp$ (see Section \ref{sectionbarabtin}).

As mentioned above, from $\mathrm{(i)}$ and (\ref{divrel}) one can derive Theorem \ref{conjmain}.
We will closely follow the ideas of Belliard.
Although we are mostly interested in the case $p=2$ or $p=3$, the method works for any prime number $p$.

\section{Elliptic units.}\label{sectionellipticunits}

For $L$ and $L'$ two $\MBZ$-lattices of $\MBC$ such that $L\subseteq L'$ and $[L':L]$ is prime to $6$, we denote by $z\mapsto\psi\left(z;L,L'\right)$ the elliptic function defined in \cite{robert90}.
For $\MFm$ a nonzero proper ideal of $\MCO_k$, and $\MFa$ a nonzero ideal of $\MCO_k$ prime to $6\MFm$, G.\,Robert proved that $\psi\left(1;\MFm,\MFa^{-1}\MFm\right) \in k(\MFm)$, where $k(\MFm)$ is the ray class field of $k$, modulo $\MFm$.
Let $S(\MFm)$ be the set of maximal ideals of $\MCO_k$ which divide $\MFm$.
Then $\psi\left(1;\MFm,\MFa^{-1}\MFm\right)$ is a unit if and only if $\vert S(\MFm)\vert = 1$.
More precisely, if we denote by $w_\MFm$ the number of roots of unity of $k$ which are congruent to $1$ modulo $\MFm$, and if we write $w_k$ for the number of roots of unity of $k$, then by \cite[Corollaire 1.3, (iv)]{robert89}, we have
\begin{equation}
\label{unitellipticisunit}
\psi\left(1;\MFm,\MFa^{-1}\MFm\right) \MCO_{k(\MFm)} = \lA \begin{array}{lll} (1) & \text{if} & 2\leq\vert S(\MFm)\vert \\ (\MFq)_{k(\MFm)}^{w_\MFm\left(N(\MFa)-1\right)/w_k} & \text{if} & S(\MFm) = \{\MFq\}, \end{array} \right.
\end{equation}
where $N(\MFa):=\#\left(\MCO_k/\MFa\right)$, and where $(\MFq)_{k(\MFm)}$ is the product of the prime ideals of $\MCO_{k(\MFm)}$ which lie above $\MFq$.
Moreover, if $\GVf_\MFm(1)$ is the Robert-Ramachandra invariant, as defined in \cite[p15]{robert73}, or in \cite[p55]{deshalit87}, we have by \cite[Corollaire 1.3, (iii)]{robert89}
\begin{equation}
\label{ellipticunitpsiandphi}
\psi\left(1;\MFm,\MFa^{-1}\MFm\right) ^{12m} = \GVf_\MFm(1)^{N(\MFa)-\left(\MFa,k(\MFm)/k\right)},
\end{equation}
where $m$ is the positive generator of $\MFm\cap\MBZ$, and $\left(\MFa,k(\MFm)/k\right)$ is the Artin automorphism of $k(\MFm)/k$ defined by $\MFa$.
If $\MFa$ is prime to $6\MFm\MFq$, then by \cite[Corollaire 1.3, (ii-1)]{robert89} we have
\begin{equation}
\label{ellipticunit}
N_{k(\MFm\MFq)/k(\MFm)} 
\left( \psi\left(1;\MFm\MFq,\MFa^{-1}\MFm\MFq\right) \right) ^{w_\MFm w_{\MFm\MFq}^{-1}} =
\lA
\begin{array}{lcl}
\psi\left(1;\MFm,\MFa^{-1}\MFm\right)^{1-\left(\MFq,k(\MFm)/k\right)^{-1}} & \mathrm{if} & \MFq\nmid\MFm,\\
\psi\left(1;\MFm,\MFa^{-1}\MFm\right) \vphantom{\left(I^J\right)^K} & \mathrm{if} & \MFq\mid\MFm.
\end{array}
\right.
\end{equation}

\begin{df}
Let $F\subseteq\MBC$ be an abelian extension of $k$, and write $\mu(F)$ for the group of roots of unity in $F$.
Let $\MFm$ be a nonzero proper ideal of $\MCO_k$.
We define the $\MBZ\left[\Gal(F/k)\right]$-submodule $\Psi(F,\MFm)$ of $F^\times$, generated by the $w_\MFm$-roots of all $N_{k(\MFm)/k(\MFm)\cap F} \left( \psi\left(1;\MFm,\MFa^{-1}\MFm\right) \right)$, where $\MFa$ is any nonzero ideal of $\MCO_k$ prime to $6\MFm$.
Also, we set $\Psi'(F,\MFm):=\MCO_F^\times\cap\Psi(F,\MFm)$.

Then, we let $C_F$ be the group generated by $\mu(F)$ and by all $\Psi'(F,\MFm)$, for any nonzero proper ideal $\MFm$ of $\MCO_k$.
\end{df}

\begin{rmq}\label{normpsooh}
If $F'\subseteq\MBC$ is a finite extension of $F$, abelian over $k$, then for all nonzero proper ideal $\MFm$ of $\MCO_k$, we have $N_{F'/F}\left(\Psi'\left(F',\MFm\right)\right) \subseteq \Psi'\left(F,\MFm\right)$ and $\Psi'\left(F,\MFm\right) \subseteq \Psi'\left(F',\MFm\right)$.
We deduce $N_{F'/F}\left(C_{F'}\right) \subseteq C_F$ and $C_F \subseteq C_{F'}$.
\end{rmq}

From (\ref{unitellipticisunit}) we deduce immediately the following remark.

\begin{rmq}\label{rmqPsirim}
If $\MFm$ is a nonzero ideal of $\MCO_k$ which is divisible by at least two distinct primes, then $\Psi'(F,\MFm)=\Psi(F,\MFm)$.
If $\MFq$ is a maximal ideal and $n\in\MBN^\ast$, then $\Psi'\left(F,\MFq^n\right)$ is generated by the products $\prod_{s\in S} N_{ k\left(\MFq^n \right) / k\left(\MFq^n \right)\cap F} \left( \psi\left(1;\MFq^n,\MFa^{-1}\MFq^n\right) \right)^{\GGa_s\GGs_s}$, where $S$ is a finite set, $\left(\MFa_s\right)_{s\in S}$ is a family of nonzero ideals of $\MCO_k$ wich are prime to $6\MFq^n$, $(\GGa_s)_{s\in S}\in\MBZ^S$ is such that $\sum_{s\in S}\GGa_s \left(N\left(\MFa_s\right)-1\right) = 0$, and $(\GGs_s)_{s\in S} \in \Gal \left( F / k \right)^S$. 
\end{rmq}

\begin{lm}\label{astrototoj}
Let $\MFm$ and $\MFg$ be two nonzero proper ideals of $\MCO_k$, such that the conductor of $F/k$ divides $\MFm$.
If $\MFg\wedge\MFm=1$, then $\Psi'\left(F,\MFg\right) \subseteq C_F\cap\MCO_{k(1)}^\times$.
Else we have $\Psi'\left(F,\MFg\right) \subseteq \Psi'\left(F,\MFg\wedge\MFm\right)$.
\end{lm}

\noindent\textsl{Proof.}
If $\MFg\wedge\MFm=1$, then $k(\MFg)\cap F\subseteq k(\MFg)\cap k(\MFm)=k(1)$.
Now assume $\MFg':=\MFg\wedge\MFm\neq1$.
There are maximal ideals $\MFq_1$, ..., $\MFq_n$ of $\MCO_k$ such that $\MFg = \MFg' \MFq_1 \cdots \MFq_n$.
By recurrence, we are reduced to the case $n=1$.
Let $\MFa$ be a nonzero ideal which is prime to $6\MFg$, and let $x$ be a $w_\MFg$-th root of $N_{k(\MFg)/k(\MFg)\cap F} \left( \psi \left( 1;\MFg,\MFa^{-1}\MFg \right) \right)$.
As above, $k(\MFg)\cap F = k(\MFg')\cap F$, and then from (\ref{ellipticunit}) there is $\GGa \in \MBZ \left[ \Gal\left( k(\MFm) / k\right) \right]$ such that 
\[x^{w_{\MFg'}} = N_{k(\MFg') / k(\MFg')\cap F} \left( N_{k(\MFg) / k(\MFg')} \left( \psi\left( 1;\MFg,\MFa^{-1}\MFg \right) \right) ^ {w_{\MFg'} / w_\MFg} \right) = N_{k(\MFg') / k(\MFg')\cap F} \left( \psi\left( 1;\MFg',\MFa^{-1}\MFg' \right) \right)^\GGa.\]
\hfill $\square$
\\

\begin{cor}\label{astrototo}
For any nonzero proper ideal $\MFm$ of $\MCO_k$, such that the conductor of $F/k$ divides $\MFm$, the group $C_F$ is generated by $\mu(F)$, by $C_F\cap \MCO^\times_{k(1)}$, and by the $\Psi'(F,\MFg)$, where $\MFg$ is a nonzero proper ideal of $\MCO_k$ which divides $\MFm$.
\end{cor}

Let $\MFf$ be the ideal of $\MCO_k$, prime to $\MFp$, such that the conductor of $K_0/k$ divides $\MFf\MFp^\infty$.
For $n\in\MBN$, and any ideal $\MFg\neq(0)$ of $\MCO_k$, we set 
\[\Psi'\left(K_n,\MFg\MFp^\infty\right) := \mathop{\cup}_{m=1}^\infty \Psi'\left(K_n,\MFg\MFp^m\right).\]
Remark that in view of (\ref{ellipticunit}), $\Psi'\left(K_n,\MFg\MFp^m\right) \subseteq \mu_{w_{\MFg\MFp^m}}(K_n)\Psi'\left(K_n,\MFg\MFp^{m+1}\right)$ for any $m\in\MBN$, 
where for any $j\in\MBN^\ast$ and any field $L$, $\mu_j(L)$ is the group of $j$-th roots of unity in $L$.
From Corollary \ref{astrototo}, we deduce the following remark.

\begin{rmq}\label{grpeellipticunitsgenerated}
The group $C_n := C_{K_n}$ is generated by $\mu(K_n)$, by $C_n\cap \MCO^\times_{k(1)}$, by $\Psi'\left(K_n,\MFp^\infty\right)$, and by the $\Psi'(K_n,\MFg)$ and the $\Psi'\left(K_n,\MFg\MFp^\infty\right)$, where $\MFg$ is a nonzero proper ideal of $\MCO_k$ which divides $\MFf$.
\end{rmq}

We define $\MCC_n:=\MBZ_p\otimes_\MBZ C_n$, $\MCC_\infty := \varprojlim \left(\MCC_n\right)$, and for any nonzero ideal $\MFg$ of $\MCO_k$, we set $\overline{\Psi}'\left(K_\infty,\MFg\MFp^\infty\right) := \varprojlim \left( \MBZ_p\otimes_\MBZ\Psi'\left(K_n,\MFg\MFp^\infty\right) \right)$.

\begin{lm}\label{grpeellipticunitsgeneratedalinfini}
The group $\MCC_\infty$ is generated by the $\overline{\Psi}'\left(K_\infty,\MFg\MFp^\infty\right)$, with $\MFg$ a nonzero ideal of $\MCO_k$ dividing $\MFf$.
\end{lm}

\noindent\textsl{Proof.}
Since $\mu\left(K_\infty\right)$ is finite, $\varprojlim\left(\mu\left(K_n\right)\right)$ is trivial.
For any proper ideal $\MFg$ dividing $\MFf$, $\Psi'\left(K_n,\MFg\right)\subseteq\MCO_{k(\MFg)}$ for all $n\in\MBN$, so one can easily check that $\varprojlim\left( \MBZ_p\otimes_\MBZ \Psi'\left(K_n,\MFg\right)\right)$ is trivial.
In the same way, $\varprojlim\left( \MBZ_p\otimes_\MBZ \left(C_n\cap \MCO^\times_{k(1)}\right) \right)$ is trivial.
Then from Remark (\ref{grpeellipticunitsgenerated}), we derive the lemma.
\hfill $\square$
\\

\begin{df}
We write $C^\SCr_F$ for the subgroup of $\MCO_F^\times$ generated by $\mu(F)$ and by the elements $N_{k(\MFm)/k(\MFm)\cap F}\left( \psi\left(1;\MFm,\MFa^{-1}\MFm\right) \right)^{\GGs-1}$, where $\MFm$ is a nonzero proper ideal of $\MCO_k$, $\MFa$ is a nonzero ideal of $\MCO_k$, prime to $6\MFm$, and $\GGs\in\Gal(F/k)$.
This is the group used for instance in \cite{rubin91} and \cite{oukhaba07}.
\end{df}

\begin{rmq}\label{sdgksgf}
It is well known that $\MCO_F^\times/C^\SCr_F$ is finite.
On the other hand, it is obvious that $C_F^\SCr\subseteq C_F$.
Hence $\MCO_F^\times / C_F$ and $C_F/C_F^\SCr$ are finite.
\end{rmq}

\begin{rmq}\label{rmqunitellipticrubin}
Using an additive notation, for $n\in\MBN$, $C_n^\SCr:=C^\SCr_{K_n}$ is generated by $\mu(K_n)$, $C_n^\SCr\cap\MCO_{k(\MFf)}^\times$, $\MBZ\left[G_n\right]_0 \Psi\left(K_n,\MFp^\infty\right)$, $\MBZ\left[G_n\right]_0 \Psi\left(K_n,\MFg\right)$ and $\MBZ\left[G_n\right]_0 \Psi\left(K_n,\MFg\MFp^\infty\right)$, where $\MFg$ is a nonzero proper ideal which divides $\MFf$, and where $\MBZ\left[G_n\right]_0$ is the augmentation ideal of $\MBZ\left[G_n\right]$.
\end{rmq}

Let $\MBQ^\alg$ be an algebraic closure of $\MBQ$.
We choose arbitrarily two embeddings $\MBQ^\alg\hookrightarrow\MBC$ and $\MBQ^\alg\hookrightarrow\MBC_p$.
For any finite abelian group $\MCG$, we denote by $\widehat{\MCG}$ the set of irreducible $\MBQ^{\alg}$ characters on $\MCG$.
For any $\chi\in\widehat{\MCG}$, we denote by $e_\chi$ the idempotent attached to $\chi$, i.e
\[e_\chi\quad=\quad\#(\MCG)^{-1}\sum_{g\in\MCG}\chi(g)g^{-1}.\]

\begin{lm}\label{quotunitellipunitrubiellip}
For $n\in\MBN$, we set $C_n^\SCr:=C^\SCr_{K_n}$ and $\MCC_n^\SCr:=\MBZ_p\otimes_\MBZ C^\SCr_n$.
For all $n\in\MBN$, the orders of $\MCC_n/\MCC^\SCr_n$ are bounded.
\end{lm}

\noindent\textsl{Proof.}
By Remark \ref{grpeellipticunitsgenerated}, we just have to show that $\left\vert C_n\cap\MCO_{k(1)}^\times/C_n^\SCr\cap\MCO_{k(1)}^\times \right\vert$ is bounded, and that the
\[X_{n,\MFg} \; := \; \MBZ_p\otimes_\MBZ\Psi'\left(K_n,\MFg\right) \; / \; \left(\MBZ_p\otimes_\MBZ\Psi'\left(K_n,\MFg\right)\right) \cap \MCC_n^\SCr \; ,\]
have bounded orders independantly of $n$ and of the nonzero proper ideals $\MFg$ of $\MCO_k$ which divide $\MFf\MFp^\infty$.
The quotient $C_n\cap\MCO_{k(1)}^\times/C_n^\SCr\cap\MCO_{k(1)}^\times$ has a bounded order because if $n$ is large enough, it is a subgroup of a quotient of $\MCO_{k(1)\cap K_\infty}^\times / C^\SCr_{k(1)\cap K_\infty}$.
In particular, there is $r\in\MBN^\ast$ an integer such that for all $n\in\MBN$, 
\begin{equation}
\label{totoOot}
\left\vert C_{k_n}\cap\MCO_{k(1)}^\times/C_{k_n}^\SCr\cap\MCO_{k(1)}^\times \right\vert\leq r.
\end{equation}

For any nonzero proper ideal $\MFg$ of $\MCO_k$ let $g\in\MBN$ be such that $\MFg\cap\MBZ=g\MBZ$.
Since $\mu_{p^\infty}\left(K_\infty\right):=\mathop{\cup}_{j=0}^\infty \mu_{p^j}\left(K_\infty\right)$ is finite, we are reduced to show that the
\[X'_{n,\MFg} \; := \; 12gw_\MFg\MBZ_p\otimes_\MBZ\Psi'\left(K_n,\MFg\right) \; / \; 12gw_\MFg \left( \left(\MBZ_p\otimes_\MBZ\Psi'\left(K_n,\MFg\right)\right) \cap \MCC_n^\SCr \right) \; \]
have bounded orders.
First, let us show that the minimal number of generators of the $\MBZ_p$-modules $X'_{n,\MFg}$ is bounded.
We notice that $X'_{n,\MFg}$ is a quotient of 
\begin{equation}
\label{fuu}
Y_{n,\MFg} \; := \; 12gw_\MFg\MBZ_p\otimes_\MBZ\Psi'\left(K_n,\MFg\right) \; / \; 12gw_\MFg I_n \left(\MBZ_p\otimes_\MBZ\Psi\left(K_n,\MFg\right)\right) \; , 
\end{equation}
where $I_n$ is the augmentation ideal of $\MBZ_p\left[G_n\right]$.
From \cite[Lemme 1.1, Chapitre IV]{tate84} and (\ref{ellipticunitpsiandphi}), we deduce that
\begin{equation}
\label{tropopoplayadurk}
Y_{n,\MFg} = \lA
\begin{array}{ll}
J_n \GVf_{n,\MFg}(1) \; / \; I_n J_n\GVf_{n,\MFg}(1) & \textrm{if} \; \MFg \; \text{is divisible by at least two distinct primes}, \\
&\\
\left(I_n\cap J_n\right)\GVf_{n,\MFg}(1) \; / \; I_n J_n \GVf_{n,\MFg}(1) & \textrm{if} \; \MFg \; \text{is a power of a maximal ideal}, \\
\end{array}
\right.
\end{equation}
where $\GVf_{n,\MFg}(1) := N_{k(\MFg) / K_n\cap k(\MFg)}\left(\GVf_{\MFg}(1)\right)$, and where $J_n$ is the annihilator of the $\MBZ_p\left[ G_n \right]$-module $\mu_{p^\infty}\left(K_n\right)$.
By (\ref{tropopoplayadurk}), it is enough to show that the minimal number of generators of the $\MBZ_p$-module $J_n/I_nJ_n$ is bounded.
It suffices to prove that the kernel of the surjections $J_n/I_nJ_n\rightarrow J_0/I_0J_0$, is bounded.
This kernel is included in $I_n\cap J_n/I_nJ_n$, which is a quotient of $I_\infty\cap J_\infty/I_\infty J_\infty$, with $I_\infty := \varprojlim \left(I_n\right)$ and $J_\infty := \varprojlim \left(J_n\right)$.
But $I_\infty\cap J_\infty/I_\infty J_\infty$ is a pseudo-nul $\GGL'$-module, annihilated by $\#\left(\mu_{p^\infty}\left(K_\infty\right)\right)$ 
and $\GGg-1$.
Hence $I_\infty\cap J_\infty/I_\infty J_\infty$ is finite.
We deduce that the minimal number of generators of the $\MBZ_p$-modules $X'_{n,\MFg}$ is bounded.

Now in order to prove the lemma, we just have to find out some $a\in\MBZ_p$ such that $a$ annihilates all the $X'_{n,\MFg}$.
Let $\GGz\in\MBC_p$ be a primitive $[K_0:k]$-th root of unity, and let $a\in\MBZ_p$ be such that $a$ is divisible by
\[r [K_0:k] \, \#\left(\mu_{p^\infty}\left(k_\infty\right)\right) \, \prod_{\substack{\chi\in\widehat{G_\MFf} \\ \chi\neq1}} \left(1-\chi\left(\GGs_\chi\right)\right)\]
in $\MBZ_p[\GGz]$, where for all nontrivial $\chi\in\widehat{G_\MFf}$, we have arbitrarily chosen $\GGs_\chi \in G_\MFf$ such that $\chi\left(\GGs_\chi\right)\neq1$.
The canonical map $X'_{n,\MFg}\rightarrow\MBZ_p[\GGz] \otimes_{\MBZ_p} X'_{n,\MFg}$ is injective, hence it is sufficient to show that $a$ annihilates $\overline{X}'_{n,\MFg} := \MBZ_p[\GGz] \otimes_{\MBZ_p} X'_{n,\MFg}$.
We have
\begin{equation}
\label{yHfrhh}
a\overline{X}'_{n,\MFg} \subseteq  r\#\left(\mu_{p^\infty}\left(k_\infty\right)\right)[K_0:k]e_1\overline{X}'_{n,\MFg} \,\oplus\, \ba{\oplus}{\substack{\chi\in\widehat{G_\MFf} \\ \chi\neq1}} \, \left(1-\GGs_\chi\right)[K_0:k]e_\chi \overline{X}'_{n,\MFg}.
\end{equation}
But obviously, $\left(1-\GGs_\chi\right)[K_0:k]e_\chi \overline{X}'_{\MFg,n}=0$ for all nontrivial $\chi\in\widehat{G_\MFf}$.
Then by (\ref{yHfrhh}) we just have to show that $r\#\left(\mu_{p^\infty}\left(k_\infty\right)\right)[K_0:k]e_1\overline{X}'_{n,\MFg}=0$.
By Remark \ref{normpsooh}, we have $\MBZ_p \otimes_\MBZ C_{k_n} \subseteq \MCC_n$ and
\begin{equation}
\label{OJHhHG}
[K_0:k]e_1\MBZ_p \otimes_\MBZ \Psi'\left(K_n,\MFg\right) \subseteq \MBZ_p \otimes_\MBZ C_{k_n}.
\end{equation}
As for (\ref{tropopoplayadurk}), for all $m\in\MBN^\ast$ we have
\begin{equation}
\label{JyYgh}
12p^mw_{\MFp^m} \Psi'\left(k_n,\MFp^m\right) / 12p^mw_{\MFp^m} I'_n \Psi\left(k_n,\MFp^m\right) \; = \; I'_n\cap J'_n\GVf'_{n,\MFp^m}(1) / I'_nJ'_n\GVf'_{n,\MFp^m}(1),
\end{equation}
where $I'_n$ is the augmentation ideal of $\MBZ\left[\GGG/\GGG_n\right]$, $J'_n$ is the annihilator of the $\MBZ\left[\GGG/\GGG_n\right]$-module $\mu_{p^\infty}\left(k_n\right)$, 
and $\GVf'_{n,\MFp^m}(1) := N_{k\left(\MFp^m\right) / k\left(\MFp^m\right)\cap k_n} \left(\GVf_{\MFp^m}(1)\right)$.
By (\ref{JyYgh}), $\#\left(\mu_{p^\infty}\left(k_\infty\right)\right)$ annihilates $\Psi'\left(k_n,\MFp^m\right) / I'_n \Psi\left(k_n,\MFp^m\right)$.
Hence by (\ref{totoOot}) and by Corollary \ref{astrototo}, we deduce that $r\#\left(\mu_{p^\infty}\left(k_\infty\right)\right)$ annihilates $C_{k_n} / C_{k_n}^\SCr$.
Then from the inclusion (\ref{OJHhHG}), we deduce the equality $r\#\left(\mu_{p^\infty}\left(k_\infty\right)\right)[K_0:k]e_1\overline{X}'_{n,\MFg}=0$, which ends the proof of the lemma.
\hfill $\square$
\\

\begin{cor}
We set $\MCC_\infty^\SCr := \varprojlim \left(\MCC_n^\SCr\right)$.
The canonical injection $\MCC^\SCr_\infty\hookrightarrow\MCC_\infty$ is a pseudo-isomorphism of $\GGL'$-modules.
\end{cor}

\noindent\textsl{Proof.}
From Lemma \ref{quotunitellipunitrubiellip}, there is $a\in\MBN$ such that $a$ annihilates $\MCC_n/\MCC^\SCr_n$ for all $n\in\MBN$.
Then $\MCC_\infty/\MCC^\SCr_\infty$ is annihilated by $a$ and by $\GGg-1$, where $\GGg$ is a topological generator of $\GGG$, hence it is pseudo-nul.
\hfill $\square$

\section{Elementary results.}\label{sectionbarabtin}

Let $F/k$ be an abelian extension of $k$.
If $[F:k]<\infty$, we write $\MCK_F$ for the pro-$p$-completion of $\prod_{\MFq\vert\MFp} F_\MFq^\times$, where for any prime ideal $\MFq$ of $\MCO_F$ above $\MFp$, $F_\MFq$ is the completion of $F$ at $\MFq$.
When $F/k$ is infinite, we define $\MCK_F$ by taking projective limits over finite sub-extensions, under the norm maps.
For all $n\in\MBN\cup\{\infty\}$, we write $\MCK_n$ for $\MCK_{K_n}$.
For any prime $\MFq$ of $K_\infty$ over $\MFp$, we define $K_{\infty,\MFq} := \mathop{\cup}_{n\in\MBN}K_{n,\MFq}$.
One can verify that $\mu_{p^\infty}\left(K_{\infty,\MFq}\right)$ is finite (otherwise see \cite[Lemma 2.1]{viguie11b}). 
Then by \cite[Theorem 25]{iwasawa73}, there is an exact sequence of $\GGL'$-modules
\begin{equation}
\label{suitexactMCKMBZGmu}
0 \rightarrow \MCK_\infty \xymatrix{\ar[r] &} \left(\GGL'\right) ^ {\left[K_0 : k\right]} \xymatrix{\ar[r] &} \ba{\oplus}{\MFq|\MFp}\mu_{p^\infty}\left(K_{\infty,\MFq}\right) \rightarrow 0.
\end{equation}
In particular $\MCK_\infty$ is a finitely generated $\GGL'$-module.

If $[F:k]<\infty$, we set $\MCE_F := \MBZ_p\otimes_{\MBZ}\MCO_F^\times$, and we write $\MCU_F$ for the pro-$p$-completion of $\prod_{\MFq\vert\MFp} \MCO_{F_\MFq}^\times$.
The injection $\MCO_F^\times \hookrightarrow \prod_{\MFq\vert\MFp} \MCO_{F_\MFq}^\times$ induces a canonical map $\MCE_F\rightarrow\MCU_F$.
The Leopoldt conjecture, which is known to be true for abelian extensions of $k$, states that this map is injective.
When $F/k$ is infinite, we define $\MCE_F$, $\MCU_F$ and a canonical injection $\MCE_F \hookrightarrow \MCU_F$, by taking projective limits over finite sub-extensions, under the norm maps.
For all $n\in\MBN\cup\{\infty\}$, we set $\MCE_n := \MCE_{K_n}$, and $\MCU_n := \MCU_{K_n}$. 
Then $\MCE_\infty$ and $\MCU_\infty$ are two submodules of $\MCK_\infty$, hence are finitely generated over $\GGL'$.

We set $B_F := \Gal\left(\MM_\MFp\left(F\right) / F\right)$, where we write $\MM_\MFp(F)$ for the maximal abelian pro-$p$-extension of $F$ which is unramified outside the primes above $\MFp$.
If $[F:k] < \infty$, let $A_F$ be the $p$-part of the class group $\Cl\left(\MCO_F\right)$.
Else, we set $A_F := \varprojlim\,A_{F'}$, where the projective limit is taken over the finite sub-extensions of $F/k$, with respect to the norm maps.
By class field theory, $A_F$ is identified to the Galois group of the maximal abelian unramified pro-$p$-extension of $F$.
For all $n\in\MBN\cup\{\infty\}$, we set $A_n := A_{K_n}$, and $B_n := B_{K_n}$. 
From \cite[end of §3.2]{iwasawa73}, we know that $B_\infty$ is a finitely generated $\GGL'$-module, hence it is also the case for $A_\infty$. 
From class field theory we have a map $\MCU_F \rightarrow B_F$, and then the following sequence of $\MBZ_p\left[\Gal(F/k)\right]$-modules, which is called the inertia sequence, is exact, 
\begin{equation}
\label{inertiasequanec}
\xymatrix{0 \ar[r] & \MCE_F \ar[r] & \MCU_F \ar[r] & B_F \ar[r] & A_F \ar[r] & 0.}
\end{equation}
If $[F:k]<\infty$, we set $\MCE'_F := \MBZ_p\otimes_\MBZ E'_F$, where $E'_F$ is the group of elements of $F$ which are units outside of the primes above $\MFp$.
When $F/k$ is infinite, we define $\MCE'_F$ by taking projective limits over finite sub-extensions, under the norm maps.
We denote by $A'_F$ the Galois group of the maximal abelian unramified pro-$p$-extension of $F$ which is totally split above $\MFp$.
For all $n\in\MBN\cup\{\infty\}$, we write $\MCE'_n$ for $\MCE'_{K_n}$, and we write $A'_n$ for $A'_{K_n}$.
The module $A'_\infty$ is a quotient of $B_\infty$, hence a finitely generated $\GGL'$-module.
As above, by the Leopoldt conjecture we have a natural injection $\MCE'_F \rightarrow \MCK_F$.
Then $\MCE'_\infty$ is finitely generated over $\GGL'$, sonce it's a submodule of $\MCK_\infty$.
From class field theory we have a map $\MCK_F \rightarrow B_F$.
Then the following sequence of $\MBZ_p\left[\left[\Gal(F/k)\right]\right]$-modules, which is called the decomposition sequence, is exact,
\begin{equation}
\label{decompossequanec}
\xymatrix{0 \ar[r] & \MCE'_F \ar[r] & \MCK_F \ar[r] & B_F \ar[r] & A'_F \ar[r] & 0.}
\end{equation}

For any $\GGL'$-module $M$, we denote by $M^{\GGG_n}$ the module of $\GGG_n$-invariants of $M$, and we denote by $M_{\GGG_n}$ the module of $\GGG_n$-coinvariants of $M$.
By definition, they are respectively the kernel and the cokernel of the multiplication by $1-\GGg_n$ on $M$, where $\GGg_n:=\GGg^{p^n}$.
We recall that if $M$ is finitely generated, then $M^{\GGG_n}$ is finite if and only if $M_{\GGG_n}$ is finite.

As in \cite[Proposition 2.2]{belliard09}, we have for all $n\in\MBN$, the exact sequence of $\GGL$-modules
\begin{equation}
\label{suitexactMCKMCKGalKGUYFYTF}
0 \xymatrix{\ar[r] &} \left(\MCK_\infty\right)_{\GGG_n} \xymatrix{\ar[r] &} \MCK_n \xymatrix{\ar[r] &} \ba{\oplus}{\MFq|\MFp}\Gal\left(K_{\infty,\MFq'}/K_{n,\MFq}\right) \rightarrow 0,
\end{equation}
where for any prime $\MFq$ of $K_\infty$ over $\MFp$, $\MFq'$ is an arbitrary prime of $K_\infty$ above $\MFq$.

Let $M_\infty:=\varprojlim\left(M_n\right)$, with $M_n$ a $\MBZ_p\left[\GGG/\GGG_n\right]$-module for every $n\in\MBN$.
For such a $\GGL'$-module, let us denote by $\Ker_n(M_\infty)$ and $\Cok_n(M_\infty)$ the kernel and cokernel of the natural map $(M_\infty)_{\GGG_n}\rightarrow M_n$, and let us denote by $\widetilde{M}_n$ the image of $M_\infty$ in $M_n$.

\begin{lm}\label{KernE'nul}
For all $n\in\MBN$, $\Ker_n(\MCE'_\infty)=0$.
\end{lm}

\noindent\textsl{Proof.}
Let $n\in\MBN$.
The $\MBZ_p$-rank of $\MCE_n$, $\MCU_n$, and $A_n$ are respectively $\left[K_n:k\right]-1$, $\left[K_n:k\right]$, and $0$.
From the inertia sequence (\ref{inertiasequanec}) applied to $K_n$, we deduce that the $\MBZ_p$-rank of $B_n$ is $1$.
We have $\left(B_\infty\right)_{\GGG_n} = \Gal\left( \MM_\MFp\left(K_n\right) / K_\infty \right)$, hence from the exact sequence below
\[\xymatrix{0 \ar[r] & \Gal\left(\MM_\MFp\left(K_n\right) / K_\infty\right) \ar[r] & B_n \ar[r] & \GGG_n \ar[r] & 0,}\]
we deduce that the $\MBZ_p$-rank of $\left(B_\infty\right)_{\GGG_n}$ is $0$.
Then $\left(B_\infty\right)_{\GGG_n}$ and $\left(B_\infty\right)^{\GGG_n}$ are finite.
By \cite[end of section 4]{greenberg78}, we know that $B_\infty$ has no nontrivial finite $\GGL'$-submodule, and we deduce
\begin{equation}\left(B_\infty\right)^{\GGG_n}=0.\label{psjdhhjg}\end{equation}
By (\ref{decompossequanec}), if we denote by $\MCK'_\infty$ the image of $\MCK_\infty$ in $B_\infty$, we have the exact sequence
\begin{equation}
\label{Olameteconouquoi}
0 \rightarrow \left(\MCE'_\infty\right)^{\GGG_n} \rightarrow \left(\MCK_\infty\right)^{\GGG_n} 
\rightarrow \left(\MCK'_\infty\right)^{\GGG_n} \rightarrow \left(\MCE'_\infty\right)_{\GGG_n} 
\rightarrow \left(\MCK_\infty\right)_{\GGG_n} \rightarrow \left(\MCK'_\infty\right)_{\GGG_n} \rightarrow 0.
\end{equation}
From (\ref{psjdhhjg}) we have $\left(\MCK'_\infty\right)^{\GGG_n}=0$.
The lemma follows from the diagram below, where the maps are injective by (\ref{Olameteconouquoi}), (\ref{suitexactMCKMCKGalKGUYFYTF}), (\ref{decompossequanec}),
\[\xymatrix{
\left(\MCE'_\infty\right)_{\GGG_n} \ar[rrr] \ar@{^{(}->}[d] & & & \MCE'_n \ar@{_{(}->}[d] \\
\left(\MCK_\infty\right)_{\GGG_n} \ar@{^{(}->}[rrr] & & & \MCK_n
}\]
\hfill $\square$
\\

Let $\MCD_\MFp$ be the decomposition group of $\MFp$ in $K_\infty/k$.
For all $n\in\MBN$, we denote by $D_n$ the subfield of $K_\infty$ fixed by $\MCD_\MFp\GGG_n$.
We set $D_\infty:=\mathop{\cup}_{n\in\MBN}D_n$ (remark that for $n$ large enough, $D_n=D_\infty$).
For any $\MBZ_p$-module $M$, we denote by $M_{\tor}$ the submodule of $M$ defined by the elements $x\in M$ such that $p^nx=0$ for large enough $n\in\MBN$.
Let us denote by $N'_{K_n/D_n}$ the composite of $N_{K_n/D_n}:\MCU_n\rightarrow\MCU_{D_n}$ with $\MCU_{D_n} \rightarrow \MCU_{D_n}/\left(\MCU_{D_n}\right)_{\tor}$ (remark that $\left(\MCU_{D_n}\right)_{\tor}$ is trivial if $p\neq2$).

For any abelian extension of local fields $L'/L$, and all $x\in L^\times$, we write $\left(x,L'/L\right)$ for the (local) norm residue symbol.
We denote by $\MCO_L^{\times,1}$ the group of units of $\MCO_L$ which are congruent to $1$ modulo the maximal ideal.

\begin{lm}\label{utildMCU}
Let $n\in\MBN$ be such that $D_\infty\subseteq K_n$, let $\MFq$ be a prime of $K_\infty$ lying above $\MFp$, and let $u \in \MCO _{K_{n,\MFq}} ^{\times,1}$.
Let us identify $\MCU_n$ to $\prod_{\MFr\vert\MFp} \MCO _{K_{n,\MFr}}^{\times,1}$, and let us write $\widetilde{\MCU}_n \left( \MFq \right)$ for the image of $\widetilde{\MCU}_n$ in $\MCO _{K_{n,\MFq}}^{\times,1}$.
Then $u\in\widetilde{\MCU}_n \left( \MFq \right)$ if and only if $N_{K_{n,\MFq} / k_\MFq} \left( u \right) \in \mathop{\cap}_{m=0}^\infty N_{k_{m,\MFq} / k_\MFq} \left( \MCO _{k_{m,\MFq}}^{\times,1} \right)$.
\end{lm}

\noindent\textsl{Proof.}
By compacity of the $\MCU_m$, we have 
\begin{equation}
\widetilde{\MCU}_n = \mathop{\cap}_{n\leq m} N_{K_m/K_n} \left(\MCU_m\right).
\label{tildeMCUint}
\end{equation}
If $u\in\widetilde{\MCU}_n \left( \MFq \right)$ then $N_{K_{n,\MFq} / k_\MFq} \left( u \right) \in \mathop{\cap}_{m=0}^\infty N_{k_{m,\MFq} / k_\MFq} \left( \MCO _{k_{m,\MFq}}^{\times,1} \right)$.
Reciprocally, assume that we have $N_{K_{n,\MFq} / k_\MFq} \left( u \right) \in \mathop{\cap}_{m=0}^\infty N_{k_{m,\MFq} / k_\MFq} \left( \MCO _{k_{m,\MFq}}^{\times,1} \right)$.
Then for $m\geq n$, we have $\left(N_{K_{n,\MFq}/k_\MFq}(u) , k_{m,\MFq} / k_\MFq\right)=1$, which implies $\left(u , K_{m,\MFq} / K_{n,\MFq}\right)_{\vert k_{m,\MFq}}=1$.
But $\Gal\left(K_m / K_n\right) \simeq \Gal\left(k_m / k_n\right)$, so the restriction map $\Gal\left( K_{m,\MFq} / K_{n,\MFq} \right) \rightarrow \Gal\left( k_{m,\MFq} / k_{n,\MFq} \right)$ is injective.
Then we deduce $\left(u , K_{m,\MFq} / K_{n,\MFq}\right)=1$, and $u\in N_{K_{m,\MFq} / K_{n,\MFq}} \left( \MCO _{K_{m,\MFq}}^{\times,1} \right)$.
Now the lemma follows from (\ref{tildeMCUint}).
\hfill $\square$
\\

\begin{lm}\label{uIUHGJJU}
Let $\MFq$ be a prime of $K_\infty$ lying above $\MFp$.
Then 
\[\mathop{\cap}_{m=0}^\infty N_{k_{m,\MFq} / k_\MFq} \left( \MCO _{k_{m,\MFq}}^{\times,1} \right) = \lA
\begin{array}{lcl}
\{1\} & \text{if} & p\neq2, \\
\mu_2 & \text{if} & p=2.
\end{array}\right.\]
\end{lm}

\noindent\textsl{Proof.}
We set $\MSN := \mathop{\cap}_{m=0}^\infty N_{k_{m,\MFq} / k_\MFq} \left( \MCO _{k_{m,\MFq}}^{\times,1} \right)$.
By local class field theory, the inertia group of $k_{\infty,\MFq} / k_\MFp$ is isomorphic to $\MCO _{k_\MFp} ^{\times,1} / \MSN$.
Since $k_{\infty,\MFq} / k_\MFp$ is infinitely ramified, we deduce that $\MSN$ is a finite subgroup of $\MCO _{k_\MFp} ^{\times,1}$.
Then $\MSN=\{1\}$ if $p\neq2$, and $\MSN\subseteq\mu_2$ if $p=2$.
Assume $p=2$.
We have $\left(-1,k_{\infty,\MFq} / k_\MFp\right)=1$ since $\GGG$ is torsion-free.
By class field theory, we deduce $-1\in\MSN$.
\hfill $\square$
\\

For all $n\in\MBN$ we define a valuation map
\[\nu_n:\mathop{\oplus}_{\MFq\vert\MFp\,\,\text{in}\,K_n} K_{n,\MFq}^\times \longrightarrow \mathop{\oplus}_{\MFq\vert\MFp\,\,\text{in}\,K_n}\MBZ\MFq, \quad x=(x_\MFq)_{\MFq\vert\MFp}\mapsto\sum_{\MFq\vert\MFp}v_\MFq\left(x_\MFq\right)\MFq,\]
where $v_\MFq$ is the normalized valuation of $K_{n,\MFq}$.
We define $\bar{\nu}_n:\MCK_n\rightarrow\mathop{\oplus}_{\MFq\vert\MFp\,\,\text{in}\,K_n}\MBZ_p\MFq$ by pro-$p$-completion, and we define $\bar{\nu}_\infty:\MCK_\infty\rightarrow\mathop{\oplus}_{\MFq\vert\MFp\,\,\text{in}\,K_\infty}\MBZ_p\MFq$ by taking projective limits.
Then we have the following exact sequence of $\GGL$-modules, for all $n\in\MBN\cup\{\infty\}$.
\begin{equation}
\label{fhffgggteh}
\xymatrix{0\ar[r] & \MCU_n \ar[r] & \MCK_n \ar[r]^-{\bar{\nu}_n} &} 
\mathop{\oplus}_{\MFq\vert\MFp\,\,\text{in}\,K_n}\MBZ_p\MFq \xymatrix{\ar[r] & 0}
\end{equation}

\begin{pr}\label{lmMCUtildeisNprimegnagna}
For all $n\in\MBN$, the $\MBZ_p$-rank of $\Ker_n\left(\MCU_\infty\right)$ and $\Cok_n\left(\MCU_\infty\right)$ is $\left[D_n:k\right]$, and we have an exact sequence of $\GGL$-modules
\begin{equation}
\label{MCUtildeisNprimegnagna}
\xymatrix{0 \ar[rr] & & \widetilde{\MCU}_n \ar[rr] & & \MCU_n \ar[rr]^-{N'_{K_n/D_n}} & & \MCU_{D_n}/\left(\MCU_{D_n}\right)_{\tor}.}
\end{equation}
\end{pr}

\noindent\textsl{Proof.}
From (\ref{suitexactMCKMBZGmu}) we know that $\left(\MCK_\infty\right)^{\GGG_n}=0$.
Then from (\ref{fhffgggteh}) we deduce the exact sequence below,
\begin{equation}
\xymatrix{0 \ar[r] &} \left(\mathop{\oplus}_{\MFq\vert\MFp\,\,\text{in}\,K_\infty}\MBZ_p\MFq\right)^{\GGG_n} \xymatrix{\ar[r] &
\left(\MCU_\infty\right)_{\GGG_n} \ar[r] & \left(\MCK_\infty\right)_{\GGG_n} \ar[r] &} 
\left(\mathop{\oplus}_{\MFq\vert\MFp\,\,\text{in}\,K_\infty}\MBZ_p\MFq\right)_{\GGG_n} \xymatrix{\ar[r] & 0.}
\label{eqliiuniyuuation}
\end{equation}
From (\ref{eqliiuniyuuation}) and (\ref{suitexactMCKMCKGalKGUYFYTF}) we deduce that the $\MBZ_p$-rank of $\Ker_n\left(\MCU_\infty\right)$ is $\left[D_n:k\right]$.

Let us show that the sequence (\ref{MCUtildeisNprimegnagna}) is exact.
We choose $m\geq n$ such that $D_\infty = D_m$.
As in \cite[Proof of Proposition 2.3]{belliard09}, the norm map $\MCU_m\rightarrow\MCU_n$ is surjective.
Then, $D_m/D_n$ being totally split at the primes above $\MFp$, we have
\begin{eqnarray}
\label{ByibYT}
\left(\MCU_{D_n}\right)_{\tor}\cap N_{K_n/D_n}\left(\MCU_n\right) & = & \left(\MCU_{D_n}\right)_{\tor}\cap N_{K_m/D_n}\left(\MCU_m\right) \nonumber\\
& = & N_{D_m/D_n} \left( \left(\MCU_{D_m}\right)_{\tor}\cap N_{K_m/D_m}\left(\MCU_m\right) \right).
\end{eqnarray}
By Lemma \ref{utildMCU} and Lemma \ref{uIUHGJJU}, we have 
\begin{equation}
\label{PIiyFyjf}
\left( \MCU_{D_m} \right)_{\tor} \cap N_{K_m/D_m}\left(\MCU_m\right) = N_{K_m/D_m} \left( \widetilde{\MCU}_m \right).
\end{equation}
From (\ref{ByibYT}) and (\ref{PIiyFyjf}), we deduce
\[ \left(\MCU_{D_n}\right)_{\tor}\cap N_{K_n/D_n}\left(\MCU_n\right) = N_{K_m/D_n} \left( \widetilde{\MCU}_m \right) = N_{K_n/D_n} \left( \widetilde{\MCU}_n \right),\]
which proves the exactness of the sequence (\ref{MCUtildeisNprimegnagna}).
Since $\MCU_{D_n}/N_{K_n/D_n}\left(\MCU_n\right)$ is a torsion $\MBZ_p$-module, we deduce that the $\MBZ_p$-rank of $\Cok_n\left(\MCU_\infty\right)$ is $\left[D_n:k\right]$.
\hfill $\square$

\section{Descent kernels and cokernels.}

For all abelian extension $F/k$, we set $\MCU_F^{(0)}:= N_{F/k}^{-1} \left( \left(\MCU_k\right)_{\tor} \right)$ (so $\MCU_F^{(0)} = \Ker\left( N_{F/k} \right)$ if $p\neq2$).
Then one can easily check that
\begin{equation}
\label{JGbbj}
\MCU^{(0)}_F/\MCE_F = \left(\MCU_F / \MCE_F\right)_{\tor} \quad \text{and} \quad \MCU^{(0)}_F/\MCC_F = \left(\MCU_F / \MCC_F\right)_{\tor}.
\end{equation}
For all $n\in\MBN$, we set $\MCU_n^{(0)}:= \MCU^{(0)}_{K_n}$, and $\MCU^{(0)}_\infty := \varprojlim\,\MCU_n^{(0)}$.
By Proposition \ref{lmMCUtildeisNprimegnagna}, we have $\widetilde{\MCU}_n\subseteq\MCU_n^{(0)}$ and therefore $\MCU^{(0)}_\infty = \MCU_\infty$.
However take care for example that for $n\in\MBN$, $\Ker_n\left(\MCU^{(0)}_\infty/\MCE_\infty\right)$ may not be equal to $\Ker_n\left(\MCU_\infty/\MCE_\infty\right)$.

\begin{lm}\label{OhFGfh}
The orders of $\Ker_n\left(\MCU^{(0)}_\infty/\MCE_\infty\right)$ and of $\Cok_n\left(\MCU^{(0)}_\infty/\MCE_\infty\right)$ are finite and bounded. 
\end{lm}

\noindent\textsl{Proof.}
By (\ref{psjdhhjg}) the inertia sequence (\ref{inertiasequanec}) gives the exact seqence below,
\begin{equation}
\label{equattgghH}
\xymatrix{
0 \ar[r] & \left(A_\infty\right)^{\GGG_n} \ar[r] & \left(\MCU_\infty / \MCE_\infty\right)_{\GGG_n} \ar[r] 
& \left(B_\infty\right)_{\GGG_n} \ar[r] & \left(A_\infty\right)_{\GGG_n} \ar[r] & 0.}
\end{equation}
Moreover, since $\MCU_n^{(0)} / \MCE_n \simeq \left(\MCU_n / \MCE_n\right)_{\tor}$ we deduce from (\ref{inertiasequanec}) the exact seqence below,
\begin{equation}
\label{equaUHHGHgghH}
\xymatrix{0 \ar[r] & \MCU_n^{(0)} / \MCE_n \ar[r] & \left(B_n\right)_{\tor} \ar[r] & A_n.}
\end{equation}
We have $\left(B_\infty\right)_{\GGG_n} \simeq \Gal\left( \MM_\MFp \left(K_n\right) / K_\infty \right) = \left(B_n\right)_{\tor}$.
Then we apply the snake lemma to the diagram below
\begin{equation}
\label{equununyt}
\xymatrix{ & 0 \ar[r] \ar[d] & \left(B_\infty\right)_{\GGG_n} \ar[r] \ar[d] & \left(B_n\right)_{\tor} \ar[r] \ar[d] & 0 \\
0 \ar[r] & \Ker_n\left(A_\infty\right) \ar[r] & \left(A_\infty\right)_{\GGG_n} \ar[r] & A_n &,}
\end{equation} 
and by (\ref{equattgghH}) and (\ref{equaUHHGHgghH}) we obtain the following exact sequence
\begin{equation}
\xymatrix{0 \ar[r] & \left(A_\infty\right)^{\GGG_n} \ar[r] & \left(\MCU_\infty / \MCE_\infty\right)_{\GGG_n} \ar[r] 
& \MCU_n^{(0)} / \MCE_n \ar[r] & \Ker_n\left(A_\infty\right) \ar[r] & 0.}
\label{equhhdeaion}
\end{equation}
As in \cite[Proof of theorem 1.4]{rubin88}, one can prove that $\left(A_\infty\right)_{\GGG_n}$ and  $\left(A_\infty\right)^{\GGG_n}$ are finite, and then from (\ref{equhhdeaion}) we deduce that the orders of $\Ker_n\left(\MCU^{(0)}_\infty/\MCE_\infty\right)$ and of $\Cok_n\left(\MCU^{(0)}_\infty/\MCE_\infty\right)$ are finite. 
Moreover $A_\infty$ is noetherian over $\GGL'$, so the sequence $\left( \Ker_n\left( \MCU^{(0)}_\infty / \MCE_\infty \right) \right)_{n\in\MBN}$ stabilizes.
For $n$ large enough (such that $K_\infty/K_n$ is totally ramified at the primes above $\MFp$), and for all $m\in\MBN$, one can easily check that $A_{n+m}\rightarrow A_n$ is surjective, and then
\begin{equation}
\Cok_n\left(A_\infty\right)=0.
\label{eioJHbfg}
\end{equation}
Moreover by a classical theorem of Iwasawa, there exists $c_1\in\MBZ$ such that for $n$ large enough, $\#\left(A_n\right) = p^{\mu\left(A_\infty\right) p^n + \GGl\left(A_\infty\right) n + c_1}$, where $\GGl\left(A_\infty\right)$ and $\mu\left(A_\infty\right)$ are respectively the $\GGl$-invariant and the $\mu$-invariant of $A_\infty$. 
Since $\left(A_\infty\right)_{\GGG_m}$ is finite for all $m\in\MBN$, it is well known, from the general theory of finitely generated $\GGL'$-modules, that there exists $c_2\in\MBZ$ such that $\#\left( \left(A_\infty\right)_{\GGG_n} \right) = p^{\mu\left(A_\infty\right) p^n + \GGl\left(A_\infty\right) n + c_2}$ for $n$ large enough.
From (\ref{eioJHbfg}), we deduce that for $n$ large enough, $\#\left( \Ker_n\left(A_\infty\right) \right) = p^{c_2-c_1}$, i.e  $\#\left( \Cok_n\left( \MCU^{(0)}_\infty/\MCE_\infty \right) \right) = p^{c_2-c_1}$ by (\ref{equhhdeaion}).
\hfill $\square$
\\

\begin{lm}\label{djhh}
The modules $\MCU_\infty / \MCC_\infty$, $\MCU_\infty / \MCE_\infty$ and $\MCE_\infty / \MCC_\infty$ are torsion over $\GGL'$, and their characteristic ideals are prime to $1-\GGg_n$, for all $n\in\MBN$.
\end{lm}

\noindent\textsl{Proof.}
It is enough to show the lemma for $\MCU_\infty / \MCC_\infty$, since $\MCU_\infty / \MCE_\infty$ is a quotient of $\MCU_\infty / \MCC_\infty$, and since $\MCE_\infty / \MCC_\infty$ is a submodule of $\MCU_\infty / \MCC_\infty$.
We refer the reader to  \cite[Proposition 3.1]{viguie11b} and \cite[Theorem 6.1]{viguie11a}. 
\hfill $\square$
\\

\begin{lm}\label{asympker}
The orders of $\Ker_n\left(\MCU^{(0)}_\infty/\MCC_\infty\right)$ and of $\Ker_n\left(\MCE_\infty/\MCC_\infty\right)$ are finite and asymptotically equivalent. 
The orders of $\Cok_n\left(\MCU^{(0)}_\infty/\MCC_\infty\right)$ and of $\Cok_n\left(\MCE_\infty/\MCC_\infty\right)$ are finite and asymptotically equivalent. 
\end{lm}

\noindent\textsl{Proof.}
From (\ref{inertiasequanec}) we deduce that $\left( \MCU^{(0)}_\infty / \MCE_\infty \right)^{\GGG_n}$ is isomorphic to a submodule of $\left(B_\infty\right)^{\GGG_n}$, hence is zero by (\ref{psjdhhjg}).
Then we have the following sequence,
\begin{equation*}
\xymatrix{ 0 \ar[r] & \left(\MCE_\infty / \MCC_\infty\right)_{\GGG_n} \ar[r] & \left(\MCU^{(0)}_\infty / \MCC_\infty\right)_{\GGG_n} \ar[r] & \left(\MCU^{(0)}_\infty / \MCE_\infty\right)_{\GGG_n} \ar[r] & 0 ,}
\end{equation*}
and we deduce the exact sequence below,
\begin{eqnarray}
& \xymatrix{0 \ar[r] & \Ker_n\left(\MCE_\infty/\MCC_\infty\right) \ar[r] & \Ker_n\left(\MCU^{(0)}_\infty/\MCC_\infty\right) \ar[r] & \Ker_n\left(\MCU^{(0)}_\infty/\MCE_\infty\right) \ar@{-}[r] & \cdots} & \nonumber\\
& \xymatrix{\cdots \ar[r] & \Cok_n\left(\MCE_\infty/\MCC_\infty\right) \ar[r] & \Cok_n\left(\MCU^{(0)}_\infty/\MCC_\infty\right) \ar[r] & \Cok_n\left(\MCU^{(0)}_\infty/\MCE_\infty\right) \ar[r] & 0.} &
\label{equationlisjhdh}
\end{eqnarray}
By Lemma \ref{djhh}, $\left(\MCU^{(0)}_\infty/\MCC_\infty\right)_{\GGG_n}$ and $\left(\MCE_\infty/\MCC_\infty\right)_{\GGG_n}$ are finite.
From Remark \ref{sdgksgf}, the group $\MCE_n/\MCC_n$ is also finite.
Then Lemma \ref{asympker} follows from (\ref{equationlisjhdh}) and Lemma \ref{OhFGfh}.
\hfill $\square$
\\

\begin{pr}\label{boundker}
The orders of $\Ker_n\left(\MCU^{(0)}_\infty/\MCC_\infty\right)$ and of $\Ker_n\left(\MCE_\infty/\MCC_\infty\right)$ are bounded.
\end{pr}

\noindent\textsl{Proof.}
By Lemma \ref{asympker}, we just have to show that the orders of $\Ker_n\left(\MCE_\infty/\MCC_\infty\right)$ are bounded.
For any $\GGL'$-module $M$, we set $\II\left(M\right) := \mathop{\cup}_{n\in\MBN} M^{\GGG_n}$.
Since the natural map $\MCE'_\infty \rightarrow \MCK_\infty$ is injective by the Leopoldt conjecture, 
we have $\left(\MCE'_\infty\right)^{\GGG_n}=0$ by (\ref{suitexactMCKMBZGmu}).
We deduce the following exact sequence,
\begin{equation}
\label{RoTuJgjg}
\xymatrix{0 \ar[r] & \left(\MCE'_\infty / \MCC_\infty\right)^{\GGG_n} \ar[r] & \left(\MCC_\infty\right)_{\GGG_n} \ar[r] & 
\left(\MCE'_\infty\right)_{\GGG_n} \ar[r] & \left( \MCE'_\infty / \MCC_\infty \right)_{\GGG_n} \ar[r] & 0.}
\end{equation}
By Lemma \ref{KernE'nul} and (\ref{RoTuJgjg}), we deduce $\Ker_n\left(\MCC_\infty\right) \simeq \left(\MCE'_\infty / \MCC_\infty\right)^{\GGG_n}$.
Since $\MCE'_\infty / \MCC_\infty$ is noetherian over $\GGL'$, we deduce that for $n\in\MBN$ large enough, we have
\begin{equation*}
\Ker_n\left(\MCC_\infty\right) \simeq \II\left(\MCE'_\infty / \MCC_\infty\right).
\end{equation*}
In the same way, for $n$ large enough we have $\Ker_n\left(\MCE_\infty\right) \simeq \II\left(\MCE'_\infty / \MCE_\infty\right)$, and then
\begin{equation}
\label{PhBHjygvgfF}
\Ker_n\left(\MCE_\infty\right) / \IM\left( \Ker_n\left(\MCC_\infty\right) \right) \simeq \II\left(\MCE'_\infty / \MCE_\infty\right) / \IM\left( \II\left(\MCE'_\infty / \MCC_\infty\right) \right).
\end{equation}

From Lemma \ref{grpeellipticunitsgeneratedalinfini} and (\ref{ellipticunit}), we deduce that $\widetilde{\MCC}_n$ is the product of the $\MBZ_p\otimes_\MBZ\Psi'\left(K_n,\MFg\MFp^\infty\right)$, where $\MFg\neq(0)$ is an ideal of $\MCO_k$ dividing $\MFf$.
From Remark \ref{grpeellipticunitsgenerated}, for all $n\in\MBN$, $\MCC_n\rightarrow\Cok_n(\MCC_\infty)$ restricts into a surjection from $\MCS_n := \MCC_n\cap\MCE_{k(\MFf)}$ onto $\Cok_n(\MCC_\infty)$.
Since $\MCE_{k(\MFf)}$ is noetherian, there is $\SCm\in\MBN$ such that for $n$ large enough, $\MCS_n = \MCS_\SCm$.
Now $\left(\MCS_\SCm\cap\widetilde{\MCC}_n\right)_{n\in\MBN}$ is an increasing sequence of submodules of $\MCS_\SCm$, which is noetherian.
Hence there exists a submodule $\MCS'$ of $\MCS_\SCm$ such that, for $n$ large enough, $\MCC_\SCm\hookrightarrow\MCC_n$ gives an isomorphism 
\begin{equation}
\label{Lhbndd}
\MCS_\SCm /\MCS' \simeq \Cok_n(\MCC_\infty).
\end{equation}
Applying the snake lemma to the diagram below,
\begin{equation*}
\xymatrix{ & \left(\MCC_\infty\right)_{\GGG_n} \ar[r] \ar[d] & \left(\MCE_\infty\right)_{\GGG_n} \ar[r] \ar[d] & \left(\MCE_\infty / \MCC_\infty\right)_{\GGG_n} \ar[r] \ar[d] & 0 \\
0 \ar[r] & \MCC_n \ar[r] & \MCE_n \ar[r] & \MCE_n / \MCC_n \ar[r] & 0,}
\end{equation*}
and since $\Ker_n\left(\MCE_\infty / \MCC_\infty\right)$ is finite, we obtain the following exact sequence,
\begin{equation}
\label{Moriarty}
\xymatrix{\Ker_n\left(\MCC_\infty\right) \ar[r] & \Ker_n\left(\MCE_\infty\right) \ar[r] & \Ker_n\left(\MCE_\infty / \MCC_\infty\right) \ar[r] & \left( \Cok_n\left(\MCC_\infty\right) \right)_{\tor}.}
\end{equation}
From (\ref{Moriarty}), (\ref{Lhbndd}), and (\ref{PhBHjygvgfF}), we deduce that for $n$ large enough, we have 
\begin{equation}
\label{Morihihihiihiarty}
\xymatrix{0 \ar[r] & \II\left(\MCE'_\infty / \MCE_\infty\right) / \IM\left( \II\left(\MCE'_\infty / \MCC_\infty\right) \right) \ar[r] & \Ker_n\left(\MCE_\infty / \MCC_\infty\right) \ar[r] & \left(\MCS_\SCm/\MCS'\right)_{\tor}.}
\end{equation}
We have a canonical injection $\MCE'_\infty / \MCE_\infty \hookrightarrow \MCU_\infty / \MCE_\infty$, and $\MCU_\infty / \MCE_\infty$ is a quotient of $\MCU_\infty / \MCC_\infty$.
By Lemma \ref{djhh}, we deduce that $\Char_{\GGL'}\left(\MCE'_\infty / \MCE_\infty\right)$ is prime to $1-\GGg_n$, hence $\left(\MCE'_\infty / \MCE_\infty\right)_{\GGG_n}$ and  $\left(\MCE'_\infty / \MCE_\infty\right)^{\GGG_n}$ are finite.
For $n$ large enough, $\II\left(\MCE'_\infty / \MCE_\infty\right) = \left(\MCE'_\infty / \MCE_\infty\right)^{\GGG_n}$, and the lemma follows from (\ref{Morihihihiihiarty}).
\hfill $\square$
\\

\begin{lm}\label{lmsudnSCm}
There exists $\SCm\in\MBN$, such that $N'_{K_n/D_n}(\MCC_n) = N'_{K_n/D_n}\left(\MCC_\SCm\cap \MCE_{k(\MFf)}\right)$ for all $n\geq\SCm$.
\end{lm}

\noindent\textsl{Proof.}
We set $|z|:=z\bar{z}$ for any $z\in\MBC$.
For any finite abelian extension $F/k$, we denote by $\ell_F$ the Dirichlet logarithm below, 
\[\ell_F: F^\times \rightarrow \MBC\left[\Gal(F/k)\right],\quad x\mapsto \sum_{\GGs\in\Gal(F/k)} \log \left\vert x^\GGs \right\vert \GGs ^{-1},\]
whose kernel verify $\Ker\left(\ell_F\right)\cap\MCO_F^\times=\mu(F)$.
From (\ref{ellipticunitpsiandphi}), we deduce 
\begin{equation}
\label{oihhjc}
\ell_{k(\MFm)} \left(\Psi'\left(k(\MFm),\MFm\right)\right) \subseteq \ba {\oplus} {\substack{ \chi\in \widehat{\Gal\left(k(\MFm)/k\right)} \\ \chi\neq1}} \MBC e_\chi\ell_{k(\MFm)} \left(\GVf_\MFm(1)\right),
\end{equation}
for any nonzero proper ideal $\MFm$ of $\MCO_k$.

Let $m\in\MBN^\ast$, and let $\MFg$ be a nonzero ideal of $\MCO_k$ which divides $\MFf$.
Let $D'$ be the subfield of $k\left(\MFg\MFp^m\right)$ fixed by the decomposition group of $\MFp$ in $k\left(\MFg\MFp^m\right)$, let $\chi\neq1$ be an irreducible complex character of $\Gal\left(D'/k\right)$, and let $\tilde{\chi}$ be the character on $\Gal \left( k\left(\MFg\MFp^m\right)/k \right)$ defined by $\chi$.
By \cite[Th\'eor\`eme 10]{robert73}, $e_{\tilde{\chi}}\ell_{k(\MFg\MFp^m)} \left(\GVf_{\MFg\MFp^m}(1)\right) = 0$.
Then from (\ref{oihhjc}) we deduce that 
\begin{equation}
\label{oihyYuuyugjc}
\ell_{k\left(\MFg\MFp^m\right)} \left(\Psi'\left(k\left(\MFg\MFp^m\right),\MFg\MFp^m\right)^{\Gal\left( k(\MFg\MFp^m) / D' \right)}\right) \subseteq \ba {\oplus} {\substack{ \chi\in \widehat{\Gal\left(D'/k\right)} \\ \chi\neq1}} \MBC e_{\tilde{\chi}}\ell_{k\left(\MFg\MFp^m\right)} \left(\GVf_{\MFg\MFp^m}(1)\right) \; = \; 0.
\end{equation}
Since $k(\MFg\MFp^m) \cap D_n \subseteq D'$, we deduce from (\ref{oihyYuuyugjc}) that $\Psi'\left(D_n,\MFg\MFp^m\right) \subseteq \mu\left(D_n\right)$.
Then $N_{K_n/D_n}\left(\widetilde{\MCC}_n\right) \subseteq \left(\MCU_{D_n}\right)_{\tor}$ by Lemma \ref{grpeellipticunitsgeneratedalinfini}. 
From Remark \ref{grpeellipticunitsgenerated}, we deduce that $N'_{K_n/D_n}(\MCC_n)$ is generated by $N'\left(\MCC_n\cap\MCE_{k(\MFf)}\right)$.
Lemma \ref{lmsudnSCm} follows because $\MCE_{k(\MFf)}$ is noetherian.
\hfill $\square$

\begin{lm}\label{fhvutmmm}
Set $\GGd:=[D_\infty:k]$.
There is $(\rho_1,\rho_2)\in\left(\MBQ_+^\ast\right)^2$ such that for all $K_n\supseteq D_\infty$,
\[\rho_1p^{(\GGd-1)n} \leq \#\left( \MCU_{D_\infty}^{(0)} / N_{K_n/D_\infty}\left(\MCU_n^{(0)}\right) \right) \leq \rho_2p^{\GGd n}.\]
\end{lm}

\noindent\textsl{Proof.}
Since $N_{K_n/D_\infty}\left(\MCU_n^{(0)}\right) = N_{K_n/D_\infty}\left(\MCU_n\right) \cap \MCU^{(0)}_{D_\infty}$, we have the exact sequence
\begin{equation}
\label{suitexacttagada}
0 \xymatrix{\ar[r] &} \frac{\MCU^{(0)}_{D_\infty}}{N_{K_n/D_\infty}\left(\MCU_n^{(0)}\right)} \xymatrix{\ar[r] &} \frac{\MCU_{D_\infty}}{N_{K_n/D_\infty}\left(\MCU_n\right)} \xymatrix{\ar[r] &} \frac{\MCU_{D_\infty}}{\MCU^{(0)}_{D_\infty}N_{K_n/D_\infty}\left(\MCU_n\right)} \longrightarrow0.
\end{equation}
Let $R_\infty$ be the maximal subextension of $K_\infty$ which is unramified over $\MFp$.
Since the quotient $\MCU_{D_\infty}^{(0)} / N_{K_n/D_\infty}\left(\MCU_n^{(0)}\right)$ is finite for all $n$ with $D_\infty\subseteq K_n$, we can assume that $R_\infty\subseteq K_n$.
From class field theory, we have an isomorphism $\frac{\MCU_{D_\infty}}{N_{K_n/D_\infty}\left(\MCU_n\right)}\simeq\prod_{\MFq|\MFp}I_\MFq$, where for all prime ideal $\MFq$ of $\MCO_{K_n}$ above $\MFp$, 
$I_\MFq$ is the $p$-part of the inertia group of $\MFq$ in $K_n/D_\infty$, which is the $p$-part of $\Gal\left(K_n/R_\infty\right)$.
The number of primes $\MFq$ in $D_\infty$ such that $\MFq|\MFp$ is $\GGd$, so if $v_1$ is the valuation at $p$ of $[K_0:k] [R_\infty:k]^{-1}$, we have
\begin{equation}
\label{equacardhdbzorn}
\#\left( \frac{\MCU_{D_\infty}}{N_{K_n/D_\infty}\left(\MCU_n\right)} \right) = p^{\GGd(v_1+n)}.
\end{equation}
The $\MBZ_p$-module $\MCU_{D_\infty}/\MCU^{(0)}_{D_\infty}$ is free of rank $1$, and so we have
\begin{equation}
\label{equacarinyfrn}
\#\left( \frac{\MCU_{D_\infty}}{\MCU^{(0)}_{D_\infty}N_{K_n/D_\infty}\left(\MCU_n\right)} \right) \leq p^{v_2+n},
\end{equation}
where $v_2$ is the valuation at $p$ of $[K_0:k] [D_\infty:k]^{-1}$.
Now Lemma \ref{fhvutmmm} follows from (\ref{suitexacttagada}), (\ref{equacardhdbzorn}), and (\ref{equacarinyfrn}).
\hfill $\square$

By Proposition \ref{lmMCUtildeisNprimegnagna}, we have $\Cok_n\left(\MCU_\infty^{(0)}\right)\simeq N'_{K_n/D_n}\left(\MCU_n^{(0)}\right)$, and the exact sequence
\begin{equation}
\label{suitexactNtildeprie}
\xymatrix{0 \ar[r] & N'_{K_n/D_n}(\MCC_n) \ar[r] & N'_{K_n/D_n}\left(\MCU_n^{(0)}\right) \ar[r] & \Cok_n\left(\MCU_\infty^{(0)}/\MCC_\infty\right) \ar[r] & 0.}
\end{equation}

\begin{pr}\label{boundcoker}
The orders of $\Cok_n\left(\MCU^{(0)}_\infty/\MCC_\infty\right)$ and of $\Cok_n\left(\MCE_\infty/\MCC_\infty\right)$ are bounded.
\end{pr}

\noindent\textsl{Proof.}
By Lemma \ref{asympker}, it is sufficient to show that the orders of $\Cok_n\left(\MCU^{(0)}_\infty/\MCC_\infty\right)$ are bounded.
Let $\SCm\in\MBN$ be as in Lemma \ref{lmsudnSCm}, large enough so that $D_\SCm=D_\infty$, and set $\MCZ:=\MCE_{k(\MFf)}\cap \MCC_\SCm$.
From Lemma \ref{lmsudnSCm} and (\ref{suitexactNtildeprie}), it is sufficient to show that the orders of $N_{K_n/D_\infty}\left(\MCU^{(0)}_n\right)/N_{K_n/D_\infty}(\MCZ)$ are bounded independantly of $n\geq\SCm$.

By (\ref{JGbbj}) and Lemma \ref{lmsudnSCm}, we see that $\MCU^{(0)}_{D_\infty}/N_{K_n/D_\infty}\left(\MCZ\right)$ is finite.
For any $x\in\MCZ$, we have $N_{K_n/D_\infty}(x)=N_{K_\SCm/D_\infty}(x)^{p^{n-\SCm}}$.
We set
\[\rho' \quad := \quad p^{(1-\GGd)\SCm} \, \#\left( \MCU^{(0)}_{D_\infty} / N_{K_\SCm/D_\infty}(\MCZ) \right) \, \#\left( \left( \MCU^{(0)}_{D_\infty} \right)_{\tor} \right),\] 
and then we have
\begin{equation}
\label{equasckepflvunt}
\#\left(\MCU^{(0)}_{D_\infty}/N_{K_n/D_\infty}(\MCZ)\right)\leq \rho' p^{(\GGd-1)n}
\end{equation}
From (\ref{equasckepflvunt}) and Lemma \ref{fhvutmmm}, we deduce
\[\#\left(N_{K_n/D_\infty}\left(\MCU^{(0)}_n\right)/N_{K_n/D_\infty}(\MCZ)\right)\leq \rho'\rho_1^{-1}.\]
\hfill $\square$

\section{Global units modulo elliptic units versus the ideal class group.}\label{Idontknow}

\begin{theo}\label{theoegalGGlmu}
The $\GGL'$-modules $\MCE_\infty/\MCC_\infty$ and $A_\infty$ share the same $\GGl$-invariant and the same $\mu$-invariant.
The $\GGL'$-modules $\MCU_\infty/\MCC_\infty$ and $B_\infty$ share the same $\GGl$-invariant and the same $\mu$-invariant.
\end{theo}

\noindent\textsl{Proof.}
By (\ref{inertiasequanec}), we have the exact sequence below,
\begin{equation*}
\xymatrix{0 \ar[r] & \MCE_\infty / \MCC_\infty \ar[r] & \MCU_\infty / \MCC_\infty \ar[r] & B_\infty \ar[r] & A_\infty \ar[r] & 0,}
\end{equation*}
from which we deduce
\begin{equation}
\label{filmdecapeetdepee}
\lA\begin{array}{lll}
\GGl\left( \MCE_\infty / \MCC_\infty \right) \GGl\left( B_\infty \right) & = & \GGl\left( \MCU_\infty / \MCC_\infty \right) \GGl\left( A_\infty \right), \\
\mu\left( \MCE_\infty / \MCC_\infty \right) \mu\left( B_\infty \right) & = & \mu\left( \MCU_\infty / \MCC_\infty \right) \mu\left( A_\infty \right). 
\end{array}\right.
\end{equation}
By (\ref{filmdecapeetdepee}), we just have to show the first assertion of the theorem.
For any two sequences $(u_n)_{n\in\MBN}$ and $(v_n)_{n\in\MBN}$ in $\MBN^\ast$, let us write $u_n\sim v_n$ if $(u_n)_{n\in\MBN}$ and $(v_n)_{n\in\MBN}$ are asymptotically equivalent.
Then by Iwasawa's theorem and \cite[Th\'eor\`eme p300, and discussion p301]{oukhaba07}, we have
\begin{equation}
\label{OuGT}
p^{\mu\left(A_\infty\right) p^n + \GGl\left(A_\infty\right) n} \, \sim \, \#\left(A_n\right) \, \sim \, \#\left(\MCE_n / \MCC_n^\SCr\right).
\end{equation}
Applying Lemma \ref{quotunitellipunitrubiellip} to (\ref{OuGT}), and then by Proposition \ref{boundker} and Proposition \ref{boundcoker}, we have
\begin{equation}
\label{camel}
p^{\mu\left(A_\infty\right) p^n + \GGl\left(A_\infty\right) n} \, \sim \, \#\left(\MCE_n / \MCC_n\right) \, \sim \, \#\left(\MCE_\infty/\MCC_\infty\right)_{\GGG_n}.
\end{equation}
By the general theory of $\GGL'$-modules, we deduce from (\ref{camel}) that
\begin{equation*}
p^{\mu\left(A_\infty\right) p^n + \GGl\left(A_\infty\right) n} \, \sim \, p^{ \mu\left( \MCE_\infty / \MCC_\infty \right) p^n + \GGl\left( \MCE_\infty / \MCC_\infty \right) n},
\end{equation*}
from which the theorem follows.
\hfill $\square$

\bibliography{biblioconjprincipdetr}
\bibliographystyle{amsplain}
\end{document}